\documentclass[12pt,reqno, oneside]{amsart}
\makeatletter
\@namedef{subjclassname@2020}{%
  \textup{2020} Mathematics Subject Classification}
\makeatother
\usepackage{amsmath,amssymb,latexsym,esint,mathrsfs,slashed, dsfont}
\usepackage{verbatim,cite,relsize,soul}
\usepackage[latin1]{inputenc}
\usepackage{microtype}
\usepackage{enumitem,graphicx, xcolor}
\usepackage{mathtools}
\mathtoolsset{showonlyrefs}
\usepackage[colorlinks=true,urlcolor=blue, citecolor=red,linkcolor=blue,
linktocpage,pdfpagelabels, bookmarksnumbered,bookmarksopen]{hyperref}
\usepackage[hyperpageref]{backref}
\usepackage[english]{babel}
\usepackage{geometry}

\numberwithin{equation}{section}
\newtheorem{theorem}{Theorem}[section]
\newtheorem{lemma}[theorem]{Lemma}

\newtheorem{proposition}[theorem]{Proposition}
\newtheorem{corollary}[theorem]{Corollary}
\theoremstyle{definition}
\newtheorem{definition}[theorem]{Definition}
\newtheorem{remark}[theorem]{Remark}

\newcommand{\R}{\mathbb R}

\title[Scattering in $\Sigma^s$ for Mass-Subcritical NLS]
{Pseudoconformal characterizations of scattering in weighted Sobolev spaces for non-scaling-invariant mass-subcritical NLS}
\author[J. Bellazzini, L. Forcella, V. Georgiev, N. Visciglia]{Jacopo Bellazzini, Luigi Forcella, Vladimir Georgiev, and Nicola Visciglia}

\address[J. Bellazzini]{Department of Mathematics, University of Pisa, Largo Bruno Pontecorvo, 5, 56127, Pisa, Italy}
\email{jacopo.bellazzini@unipi.it}

\address[L. Forcella]{Dipartimento di Matematica, Universit\`a Degli Studi di Pisa, Largo Bruno Pontecorvo, 5, 56127, Pisa, Italy}
\email{luigi.forcella@unipi.it}

\address[V. Georgiev]{Dipartimento di Matematica, Universit\`a Degli Studi di Pisa, Largo Bruno Pontecorvo, 5, 56127, Pisa, Italy, and Faculty of Science and Engineering, Waseda University, 3-4-1, Okubo, Shinjuku-ku, Tokyo 169-8555, Japan, and Institute of Mathematics and Informatics,  Bulgarian Academy of Sciences,  Acad. Georgi Bonchev Str., Block 8, 1113 Sofia, Bulgaria}
\email{vladimir.simeonov.gueorguiev@unipi.it}

\address[N. Visciglia]{Dipartimento di Matematica, Universit\`a Degli Studi di Pisa, Largo Bruno Pontecorvo, 5, 56127, Pisa, Italy}
\email{nicola.visciglia@unipi.it}



\begin{document}
	
\begin{abstract} We consider mass-subcritical nonlinear Schr\"odinger equations with finitely many power-type nonlinearities and initial data in weighted Sobolev spaces. In the short-range regime, we prove a pseudoconformal characterization of scattering in the weighted energy space $\Sigma$. We also establish a fractional counterpart in $\Sigma^s$, $s\in(0,1)$, formulated in terms of the boundedness of the order-$s$ pseudoconformal operator. Our results include equations with combined nonlinearities, for which no distinguished scaling symmetry is available. As an application of the fractional theory, we prove small-data scattering in $\Sigma^s$, $s\in(0,1]$. The arguments rely on the lens transformation, harmonic-oscillator Strichartz estimates, and fractional harmonic Sobolev estimates.
\end{abstract}	
	
\maketitle

\section{Introduction and main results}
In this paper, we consider the Cauchy problem for the  nonlinear Schr\"odinger equation (NLS) 
\begin{equation}\label{NLS0}
    \begin{cases}
    i \partial_t \psi +\Delta\psi =  \sum_{j=1}^N\mu_j|\psi|^{p_j-1}\psi\\
     \psi(0)= \psi_0  \in \Sigma
    \end{cases},
\end{equation}
where $\psi=\psi(t,x)$, $\psi:\mathbb{R}\times \R^d \to \mathbb{C}$,  $\Delta$ is the standard Laplace operator with respect to the spatial variable $x$,  and for any $j\in\{1,\dots,N\}$, $\mu_j\in \mathbb R$ and 
\begin{equation}\label{range-p}
    1< p_j<1+\frac4d.
    \end{equation}
The space $\Sigma$ is defined as the subspace of $H^1(\mathbb R^d)$ with finite variance:
\[
\Sigma=\{f\in H^1(\mathbb R^d)\, \hbox{ s.t. } \, |x|f \in L^2(\mathbb R^d) \},
\]
endowed with norm $\|f\|_{\Sigma}^2=\|f\|_{H^1(\mathbb R^d)}^2+\||x|f\|_{L^2(\mathbb R^d)}^2$. \medskip 

In the case of a single nonlinear term, i.e., $N=1$, equation 
\eqref{NLS0} reduces to 
\begin{equation}\label{NLS-s}
    \begin{cases}
    i \partial_t \psi +\Delta\psi = \mu|\psi|^{p-1}\psi\\
     \psi(0)= \psi_0  \in \Sigma
    \end{cases},
\end{equation}
which  is invariant under the scaling   
\begin{equation}\label{NLS-s-scaling}
  \psi_\lambda(t,x)
  :=\lambda^{\frac{2}{p-1}}\psi(\lambda^2t,\lambda x),
  \qquad \lambda>0,
\end{equation}
and the  $\dot H^{s_c}(\mathbb R^d)$-norm with  $s_c=\frac{d}{2}-\frac{2}{p-1}$ is  invariant under the scaling \eqref{NLS-s-scaling}. Note that $s_c=0$, corresponding to the critical space $L^2(\mathbb R^d)$,  if and only if $p=1+\frac4d$. Thus, regardless of the number of nonlinear terms,  \eqref{range-p} means we are considering  a mass-subcritical regime. 
It is well known that, for any $\psi_0\in\Sigma$, there exists  a unique solution  $\psi(t)\in C(I_{\max};\Sigma)$ to \eqref{NLS0}, and due to the mass-subcriticality assumption, the solution is global, i.e., $I_{\max}=\mathbb R$, see \cite{C03}, and by the blowup alternative $\sup_{t\in\R}\|\psi(t)\|_{H^1(\R^d)}\lesssim1$. Furthermore, the following quantities are preserved along the flow:  the mass, 
\begin{equation*}\label{eq:mass}
    M(\psi(t)) = \int_{\mathbb{R}^d} |\psi(t,x)|^2 dx,
\end{equation*}
and  the energy
\begin{equation*}\label{eq:en}
     E(\psi(t)) = \frac12\int_{\mathbb{R}^d} |\nabla \psi(t,x)|^2 dx 
   + \sum_{j=1}^N\frac{\mu_j}{p_j+1} \int_{\mathbb{R}^d} |\psi(t,x)|^{p_j+1} dx.
\end{equation*}
Once global solutions are known to exist, it is natural to investigate their long-time behaviour. Roughly speaking, the notion of scattering means that a nonlinear flow behaves as a free wave asymptotically in time. More precisely, we use the following definition of scattering. In what follows, $e^{it\Delta}$ denotes the free Schr\"odinger propagator.   The scattering theory for nonlinear Schr\"odinger equations has a long history. Foundational results on the Cauchy and scattering theories were obtained by Ginibre and Velo in \cite{GV79} and subsequent works. Wave operators in weighted Sobolev spaces were studied by Cazenave and Weissler \cite{CW92}.  Ginibre, Ozawa, and Velo \cite{GOV94} subsequently developed a fractional weighted-space approach based on Besov norms associated with the pseudoconformal vector field. 
\smallskip 

The purpose of the present paper is different from the construction of wave operators. We give an intrinsic characterization of scattering for a given global solution in the weighted energy space $\Sigma$. More precisely, we prove that scattering in $\Sigma$ is equivalent to a uniform first-order pseudoconformal bound along the nonlinear evolution. \smallskip

Let us give the rigorous definition of scattering. 
\begin{definition} Let $\psi\in C(\mathbb R;Y)$ be a global solution to \eqref{NLS0} on a Banach space $Y=Y(\mathbb R^d)$. We say that $\psi(t)$ scatters in $Y$ forward (respectively, backward) in time, provided there exists a scattering state $\psi^+\in Y$ (resp. $\psi^-\in Y$) such that 
\begin{equation}\label{def-scattering}
    \lim_{t\to+ \infty} \|U(-t)\psi(t)-\psi^+\|_Y=0,  \quad (\text{resp.} \lim_{t\to- \infty} \|U(-t)\psi(t)-\psi^-\|_Y=0)
\end{equation}
\end{definition}
\noindent where $U(t)=e^{it\Delta}$ is the free Schr\"odinger operator. If $Y=L^2(\R^d)$ or $H^1(\R^d)$, the unitary property of $U(t)$ implies that  \eqref{def-scattering} is equivalent to the convergence of the nonlinear flow to a free wave, i.e., $    \lim_{t\to\pm \infty} \|\psi(t)-U(t)\psi^\pm\|_Y=0$. In $\Sigma$, however, $U(t)$ is not unitary with respect to the standard $\Sigma$-norm, and its operator norm grows with $|t|$. Consequently, the two formulations are not immediately equivalent. See \cite{Beg} for a positive result in this direction. Scattering and wave operators in weighted Sobolev spaces were studied in several earlier works, including \cite{T85,CW92,NONO02}. In particular, pseudoconformal methods and the reduction to a nonautonomous equation on a finite time interval play an important role in \cite{CW92}. We refer also to \cite{NONO02} for a unified treatment of scattering methods in Sobolev and weighted spaces.\medskip

The main purpose of this paper is to identify a direct dynamical criterion for
scattering in the weighted energy space $\Sigma$. The criterion is expressed
through the pseudoconformal vector field 
\[
J(t)=x+2it\nabla,
\]
or, equivalently, through the gradient of the quadratic-phase-modulated solution. This viewpoint
is especially useful for equations with several nonlinear powers, since in
that setting there is generally no common scaling symmetry and hence no single
scale-invariant weighted space in which to formulate the scattering problem.

\subsection{Main results} 
In what follows, we consider the positive time direction. Hence, scattering means \emph{scattering forward in time}. For negative time, the analysis is the same. Our first main result is the following characterization of scattering in $\Sigma$.
\begin{theorem}\label{thm:char}Let $1+\max\left\{\frac2d,\frac{4}{d+2}\right\}<p_j<1+\frac4d$, $j\in\{1,\dots, N\}$. Let $\psi\in C( [0,\infty);\Sigma)$ be a global solution to \eqref{NLS0}. Then the following characterization of scattering in $\Sigma$ holds: 
\[
\sup_{t\geq 0} \|J(t)\psi(t)\|_{L^2(\R^d)}<\infty \iff \psi(t) \text{ scatters in } \Sigma. 
\]
\end{theorem}
Theorem~\ref{thm:char} separates the scattering problem into two complementary
directions. The necessity of the pseudoconformal bound follows directly from
the intertwining identity for the free Schr\"odinger group. The converse is based on the lens transformation of \eqref{NLS0}, which converts the long-time scattering problem for the free
Schr\"odinger flow into a finite-time problem for the harmonic flow. Moreover, working in $\Sigma$, the harmonic flow provides the natural framework to prove scattering in the $\Sigma$-topology.    
 \begin{remark}  The short-range assumption $p_j>1+\frac2d$ is natural in view of the failure of ordinary $L^2$-scattering for nontrivial solutions in the long-range regime $1<p\leq1+\frac2d$; see \cite{Barab,C03}. Scattering and wave operators for mass-subcritical NLS in weighted spaces have been studied by several complementary methods. By exploiting the pseudoconformal transformation, $L^2$-scattering was obtained by Tsutsumi and Yajima \cite{TY84}; see also \cite{T85}. More recently, the corresponding scattering problem in the $H^1$-topology was treated in \cite{BGTV}. Wave operators in $\Sigma$ and related weighted spaces were studied in \cite{CW92}.   Ginibre, Ozawa, and Velo \cite{GOV94} subsequently developed a fractional weighted theory based on Besov spaces associated with the pseudoconformal vector field, proving the existence of wave operators and small-data asymptotic completeness in a broad short-range regime.\smallskip
    
    \noindent The result above has a different character: it does not construct a solution from a prescribed asymptotic state, but characterizes scattering for a given global solution through the uniform first-order bound $\sup_{t\geq0}\|J(t)\psi(t)\|_2<\infty$. The proof is based on the lens transformation and endpoint continuity of the corresponding harmonic evolution. \end{remark}

  We now give a fractional generalization of \ref{thm:char}. For $s\in\mathbb R$, let us denote the fractional $\Sigma^s$-spaces 
\begin{equation}\label{def:sigma-frac}
\Sigma^{s}:=\{f \text{ s.t. } |\nabla|^sf\in L^2(\R^d) \text{ and } \langle x\rangle^sf\in L^2(\R^d)\}
\end{equation}
equipped with the norm
\begin{equation}\label{def:s-norm}
\|f\|_{\Sigma^s}^2=\|f\|_{L^2}^2+\||\nabla|^s f\|_{L^2}^2+\||x|^sf\|_{L^2}^2. 
\end{equation}
Let us introduce the order-$s$ pseudoconformal operator
 \begin{equation}\label{def-Js} 
 |J(t)|^s=U(t)|x|^sU(-t)
 \end{equation}
 Recall that $J(t)=U(t)xU(-t)$, so \eqref{def-Js} is a natural fractional extension of $J(t)$. 
 We have the following.  
 \begin{theorem}\label{thm:char-frac} Let $s\in(0,1)$, and suppose that $ 1+\max\left\{ \frac{2}{d}, \frac{4}{d+2s} \right\} < p_j < 1+\frac{4}{d}$, $ j\in\{1,\dots,N\}$.  Let $ \psi\in C([0,\infty);\Sigma^s) $ be a global solution to \eqref{NLS0} with \begin{equation}\label{cond-bound}\sup_{t\geq0}\|\psi(t)\|_{H^s(\mathbb R^d)}<\infty.\end{equation} Then the following characterization of forward scattering in $\Sigma^s$ holds: \[ \sup_{t\geq0}  \||J(t)|^s\psi(t) \|_{L^2(\mathbb R^d)} <\infty  \iff  \psi(t) \text{ scatters forward in} \Sigma^s.
 \]
\end{theorem}
\begin{remark} The characterization in Theorem~\ref{thm:char-frac} is conditional, since the uniform bound \eqref{cond-bound} is an additional hypothesis and does not follow, in general, merely from the global existence of $\psi\in C([0,\infty);\Sigma^s)$. This is the main difference from the corresponding characterization in $\Sigma$, where the required uniform $H^1$-bound is provided by the global theory for the mass-subcritical equation. In the fractional setting, the uniform $H^s$-bound is needed, together with the bound on $\||J(t)|^s\psi(t)\|_{L^2(\mathbb R^d)}$, to obtain the uniform estimate \[ \sup_{0\leq\tau<\pi/4}\|\varphi(\tau)\|_{\Sigma^s}<\infty \] for the lens-transformed solution.\end{remark}
\begin{remark} For the single-power equation \eqref{NLS-s}, scattering in critical weighted spaces $\mathcal F\dot H^{s_c}(\mathbb R^d)$ was studied in \cite{KMMV17}; the present result instead concerns scattering in $\Sigma^s$ and also applies to the non-scaling invariant equation \eqref{NLS0}. \end{remark}

Our last  aim is to  develop a small-data theory in
 low regularity weighted spaces.  We prove the following. 
\begin{theorem}\label{thm:small-data-frac}
Fix $s\in(0,1]$. Assume that  
$1+ \max\left\{\frac2d,\frac4{d+2s}\right\}<p_j<1+\frac4d$, $j\in\{1,\dots, N\}$.
Then there exists $\varepsilon=\varepsilon(d,s,p_1,\dots,p_N, \mu_1,\dots, \mu_N)>0$ such that for every $\psi_0\in\Sigma^s$ satisfying
$
  \|\psi_0\|_{\Sigma^s}\leq\varepsilon$, the unique global 
 solution $
  \psi\in C(\R;\Sigma^s)
$
scatters forward in time in the sense of \eqref{def-scattering} in $\Sigma^s$.
Thus, the forward scattering map is defined in a neighborhood
of the origin in $\Sigma^s$.
\end{theorem}
We comment about the above result.
\begin{remark} 
\textup{(i)} For $s\neq1$, unlike the characterization result of Theorem \ref{thm:char-frac}, the  small-data scattering theorem above is unconditional: the uniform $\Sigma^s$-bound for the lens-transformed solution is obtained directly from the smallness of the initial datum through a bootstrap argument.   
\smallskip

\noindent\textup{(ii)} We give here  a direct and self-contained proof in the spaces $\Sigma^s$, based on the lens transformation and fractional harmonic Sobolev estimates. \smallskip
 
\noindent \textup{(iii)} Note that the conclusion of  Theorem \ref{thm:small-data-frac} applies when $s=1$, and for  this latter value the space $\Sigma^1$ is equivalent to the usual finite-variance subspace of $H^1(\mathbb R^d)$.     Hence, for $s=1$, we recover the classical small-data scattering theory in $\Sigma$; see, for instance, \cite{T85,CW92,GOV94,NONO02,C03}.  
\end{remark}
\begin{remark} The threshold $p=1+\frac{4}{d+2s}$ is natural from the point of view of scaling. Indeed, below this threshold one can construct arbitrarily small standing-wave data in $\Sigma^s$, as the following argument shows. Consider $N=1$, $p=p_1<\frac{d+4+2s}{d+2s}$, and let 
 $Q$ be a positive radial solution to $-\Delta Q+Q-Q^{p}=0$. The rescaled profiles $Q_\lambda(x)=\lambda^\frac {2}{p-1}Q(\lambda x)$ generate standing-wave solutions with a suitably rescaled frequency.  Moreover, the norm 
\begin{equation*}\label{eq:scaling-sigma}
\begin{aligned}
\|\lambda^\frac 2{p-1}Q(\lambda \cdot)\|_{\Sigma^s}^2&= 
\lambda^{\frac{4}{p-1}+2s-d} \||\nabla|^s Q\|_{L^2(\R^d)}^2\\
&+ \lambda^{\frac{4}{p-1}-2s-d}\||\cdot|^s Q\|_{L^2(\R^d)}^2+\lambda^{\frac{4}{p-1}-d}\|Q\|_{L^2(\mathbb R^d)}^2,
\end{aligned}
\end{equation*}
 can be  arbitrarily small provided that $\lambda\to0$ for 
 $p<\frac{d+4+2s}{d+2s}$. The endpoint itself requires a separate analysis; related endpoint small-data results in fractional weighted spaces were obtained in \cite{NONO02}.
\end{remark}
As a consequence of Theorem  \ref{thm:small-data-frac}, we obtain the following scattering result for the scaling-invariant equation \eqref{NLS-s}, regardless of whether the nonlinearity is defocusing or focusing.
\begin{corollary}\label{coro-sd}
Let $s\in(0,1]$ and assume that  
$1+ \max\left\{\frac2d,\frac4{d+2s}\right\}<p<1+\frac4d$. Set $\theta=\frac{4}{p-1}-d$. Let $\psi_0\in\Sigma^s$. There exists a $\varepsilon>0$ sufficiently small such that, if \[
\||\nabla|^s\psi_0\|_2^{\frac12-\frac{\theta}{4s}}\||x|^s\psi_0\|_2^{\frac12+\frac{\theta}{4s}}<\varepsilon,
\]
then the solution to \eqref{NLS-s} with initial datum $\psi_0$ scatters forward $\Sigma^s$. 
\end{corollary}
\begin{remark}
    Corollary \ref{coro-sd} also implies the existence of scattering solutions with arbitrarily large $\Sigma^s$-norm in the focusing case.
\end{remark}

\section{Preliminary results and tools}
Since now on, we will omit the domain dependence when we are dealing with space functions defined on the whole space $\mathbb R^d$. Moreover, we will shorten the Lebesgue norm notation in the following way: $\|f\|_{p}=\|f\|_{L^p(\mathbb R^d)}$. 
\subsection{Fractional spaces and Strichartz estimates} 
Given the harmonic oscillator operator $\mathcal H=-\Delta+|x|^2$,  we introduce the space
\begin{equation*}\label{def:frac-space}
\mathcal W^{s,q}:=\{f\in L^q  \text{ s.t. } \mathcal H^{s/2}f\in L^q\}, \quad s\in\mathbb R,
\end{equation*}
endowed with the norm $\|f\|_{\mathcal W^{s,q}}^q:=\|f\|_q^q+\|{\mathcal H}^{s/2}f\|_q^q$. 
By  \cite[Theorem 1.4]{Dziu},  for $q>1$, we have the  norm equivalence 
\begin{equation}\label{norm-equiv2}
\|\mathcal H^{s/2}f\|_{q}^q\sim \||\nabla|^s f\|^q_{q}+\|\langle x\rangle^sf\|^q_q.
\end{equation}
In particular, for $\Sigma^s$ as defined in \eqref{def:sigma-frac} and \eqref{def:s-norm}, 
by  \eqref{norm-equiv2}, 
\[\Sigma^s\sim \mathcal W^{s,2}.\]
Recall the definition of the Strichartz admissible pair. A pair $(\gamma,\rho)$ is Schr\"odinger admissible if
\begin{equation}\label{eq:admissibility}
  \frac2\gamma=d\left(\frac12-\frac1\rho\right),
  \qquad
    (\gamma, \rho, d)\neq(2,\infty,2),
\end{equation}
where 
\begin{equation*}
2\leq\rho\begin{cases}
\leq \infty \quad d=1\\
<\infty \quad d=2\\
\leq 2^*:= \frac{2d}{d-2} \quad d\geq 3
\end{cases}.
\end{equation*}
For a solution $\psi$ to 
\[
\begin{cases}
i\partial_t\psi+\Delta \psi=F\\
\psi(0,x)=\psi_0
\end{cases}
\]
and $I\ni 0$, we have 
\begin{equation*}\label{eq:st-est}
\|\psi\|_{L^{\gamma}_t(I; W^{1,\rho}_x)}\lesssim  \|\psi_0\|_{H^1}+\|F\|_{L^{\tilde\gamma'}_t(I;W^{1,\tilde\rho'}_x)}
\end{equation*}
for any $(\gamma, \rho)$ and $(\tilde\gamma,\tilde\rho)$ admissible. \smallskip 

\noindent Consider now a   solution $\varphi$ to 
\[
\begin{cases}
i\partial_t\varphi-\mathcal H \varphi=F\\
\varphi(0,x)=\varphi_0
\end{cases}.
\]
The homogeneous and inhomogeneous Strichartz estimates associated with $e^{it\mathcal H}$  extend to $\mathcal W^{s,q}$. Eventually, by \eqref{norm-equiv2}, we have 
\begin{equation*}\label{eq:st-est2}
\|\varphi\|_{L^{\gamma}_t(I; \mathcal W^{s,\rho}_x)}\lesssim  \|\varphi_0\|_{\Sigma^s}+\|F\|_{L^{\tilde\gamma'}_t(I;\mathcal W^{s,\tilde\rho'}_x)}
\end{equation*}
for any $(\gamma, \rho)$ and $(\tilde\gamma,\tilde\rho)$ satisfying \eqref{eq:admissibility}.

\subsection{The Lens transformation} 

Given a time-independent function $f(x)$, we define, see \cite{BT24, BTZ, Carles09} the lens transform operator $\mathcal L_t$ as the scaling and modulation operator  which acts on $f$ as follows:
\begin{equation*}\label{lens-def}
\mathcal L_t f(x)=(\cos(2t))^{-\frac d2} f\left(\frac{x}{\cos(2t)}\right)e^{-i\frac{|x|^2\tan(2t)}{2}}, \quad t\in\left(-\frac\pi4,\frac\pi4\right), \quad x\in \mathbb R^d.
\end{equation*}
This  transformation can be viewed as a pseudoconformal compactification of the Schr\"odinger flow; see, for instance, \cite{Carles09, BGTV, Tao09}.\smallskip

\noindent  We have the following identity:
\begin{equation}\label{lens-ide-porp}
e^{-it_s\mathcal{H}}f=\mathcal L_{t_s}(e^{is\Delta}f), \quad \hbox{ with } \quad t_s:=\frac{\arctan(2s)}{2}.
\end{equation}
The lens transformation maps
infinite Schr\"odinger time to the finite interval
$(-\pi/4,\pi/4)$ and replaces the free Laplacian by the harmonic oscillator. Specifically,  the lens transformation  maps a solution   $\psi$ to \eqref{NLS0}, to a solution 
\begin{equation}\label{eq:lens-varphi}
  \varphi(t_s,x):=\mathcal L_{t_s}(\psi(s,\cdot))(x)
\end{equation} to the following Cauchy problem 
\begin{equation}\label{NLS-lens}
    \begin{cases}
    i \partial_t \varphi -\mathcal{H}\varphi =  \sum_{j=1}^N\mu_j( \cos(2t))^{-\delta(p_j)}|\varphi|^{p_j-1}\varphi,\quad (t,x)\in\left(-\frac\pi4,\frac\pi4\right)\times \mathbb R^d\\
     \varphi(0)= \psi_0 
    \end{cases},
\end{equation}
where   
$\delta(p_j)=\frac{4-d(p_j-1)}{2}$. 
By denoting by 
\[
g(t)=\sum_{j=1}^N\mu_j( \cos(2t))^{-\delta(p_j)}|\varphi|^{p_j-1}\varphi\]
the nonlinearity in \eqref{NLS-lens}, the integral formulation of the solution $\varphi$ to \eqref{NLS-lens}
is given by the usual Duhamel's formula, i.e.,  
\begin{equation}\label{eq-duha}
    \varphi(t) = e^{-it\mathcal{H}} \varphi_0 -i \int_0^t e^{-i(t-s)\mathcal{H}} g(s) ds. 
\end{equation}
The next result provides a characterization of scattering in $\Sigma^s$ for solutions $\psi(t)$ to \eqref{NLS0}: it is enough to show that the transformed function $\varphi$ admits a continuous extension in $\Sigma^s$ up to $t=\frac\pi4$. 
Accordingly, scattering is reduced to uniform control of the transformed
solution up to the endpoint of this finite interval, as the next proposition shows. 
\begin{proposition}[see \cite{BGTV}]\label{prop:equiv-sca} Let $\psi_0\in \Sigma^s$, $s\in(0, 1]$, and $p_j$ satisfy \eqref{range-p}. Then the forward $\Sigma^s$-scattering of $\psi(t)$, as in definition \eqref{def-scattering}, is equivalent to the existence of a profile $\varphi^+\in \Sigma^s$ such that 
\begin{equation}\label{def-sca-lens}
\lim_{\tau\to{\frac{\pi}{4}}^-}\|\varphi(\tau)-\varphi^+\|_{\Sigma^s}=0,
\end{equation}
where $\varphi(\cdot)$ is the solution to \eqref{NLS-lens} with corresponding initial datum $\psi_0$.
\end{proposition}
\begin{proof}
    The proof is analogous to that of \cite[Proposition 6.1]{BGTV}. For the sake of completeness, we report the proof here adapted to the fractional setting. By \eqref{lens-ide-porp}, we have that $e^{-i\sigma\Delta}=e^{it_\sigma\mathcal H}\mathcal L_{t_\sigma}$. Therefore, by \eqref{eq:lens-varphi}
    \begin{equation}\label{eq-ide-fund}
  e^{-i\sigma\Delta}\psi(\sigma)= e^{it_\sigma\mathcal H} \left( \mathcal L_{t_\sigma}\psi(\sigma) \right) =e^{it_\sigma\mathcal H}\varphi(t_\sigma).     
    \end{equation}
    Note that by the spectral calculus $\mathcal H^{s/2}$ commutes with $e^{it\mathcal H}$, hence the propagator $e^{it\mathcal H}$ is unitary  on the harmonic Sobolev space $\mathcal W^{s,2}$, whose norm is equivalent to the $\Sigma^s$-norm.  The conclusion follows by the isometry property and by  $\lim_{\sigma\to\infty}t_\sigma=\frac\pi4$.
\end{proof}

An important ingredient towards the proof of Theorem \ref{thm:char} and Theorem \ref{thm:char-frac} is the following.

\begin{proposition}\label{prop:control-sigma-frac} Let $0<s\leq1$, and let $\psi\in C([0,\infty);\Sigma^s)$ be a global solution to \eqref{NLS0}. Assume that $ \sup_{t\geq0}\||J(t)|^s\psi(t)\|_2<\infty$. Moreover, for $s\neq1$  also assume that $\sup_{t\geq0}\|\psi(t)\|_{H^s}<\infty$. Then the corresponding lens-transformed solution satisfies \[ \sup_{0\leq\tau<\pi/4} \|\varphi(\tau)\|_{\Sigma^s}<\infty. \] \end{proposition} \begin{proof} 
By \eqref{eq-ide-fund}, with $\tau=t_\sigma$,  the isometry property of $e^{i\tau \mathcal H}$ in $\mathcal W^{s,2}\sim \Sigma^s$,  the isometry property of $e^{it\Delta}$ in $L^2$, and \eqref{def-Js}, we have
\begin{equation}
    \begin{aligned}
 \|\varphi(\tau)\|_{\Sigma^s}&= \|e^{i\tau \mathcal H} \varphi(\tau)\|_{\Sigma^s}\sim \|e^{i\tau \mathcal H} \varphi(\tau)\|_{2}+ \|e^{i\tau \mathcal H} \varphi(\tau)\|_{\dot H^s} + \||x|^s e^{i\tau \mathcal H} \varphi(\tau)\|_{2}\\
 &=\| e^{-i \sigma\Delta} \psi(\sigma)\|_{2}+ \|e^{-i \sigma \Delta} \psi(\sigma)\|_{\dot H^s} + \||x|^s e^{-i \sigma\Delta} \psi(\sigma)\|_{2}\\
 &=\|\psi(\sigma)\|_{2}+ \|\psi(\sigma)\|_{\dot H^s}+ \||J(\sigma)|^s \psi(\sigma)\|_{2},
    \end{aligned}
    \end{equation}
    and the result follows by the  hypothesis. 
\end{proof}

\subsection{Fractional uncertainty principle}
In this section, we give a short and self-contained proof of the fractional uncertainty principle, see \eqref{frac_unc-2}. Althought this estimate is classical and well-known, we provide here a direct proof based only on the Poincar\'e inequality. 
First, we recall the definition of the seminorm $\dot H^s(B)$, for $s\in(0,1)$, where $B$ is the unit ball in $\mathbb R^d$:
\[
\|f\|^2_{\dot H^s(B)}=\int_B\int_B\frac{|f(x)-f(y)|^2}{|x-y|^{d+2s}}dxdy.
\]
\begin{proposition}\label{prop:u-p} Let $d\geq 1$, $s\in(0,1]$.   There exists a constant $C=C(d,s)>0$ such that \begin{equation}\label{frac_unc} \|u\|_{L^2(\mathbb R^d)}^2 \leq C\left( \|u\|_{\dot H^s(B)}^2 + \||x|^s u\|_{L^2(\mathbb R^d)}^2 \right).\end{equation} 
In particular, the following also holds:
 \begin{equation}\label{frac_unc-2} \|u\|_{2}^2 \leq C\||\nabla|^s u\|_{2}\||x|^su\|_{2}.\end{equation} 
 \end{proposition} \begin{proof} 
In the case $s=1$, \eqref{frac_unc-2} is the classical Heisenberg inequality. For $s\in(0,1)$,  set $ u_B:=\frac{1}{|B|}\int_B u(x)\,dx$.
By the fractional Poincar\'e inequality, see \cite{Hurri}, there exists $C_{d,s}>0$ such that \[ \|u-u_B\|_{L^2(B)} \leq C_{d,s}\|u\|_{\dot H^s(B)}. \]  Consider the annulus \[ A:=\left\{x\in\mathbb R^d:\frac12<|x|<1\right\}. \] 
Since $A\subset B$, we have \begin{equation}\label{eq:1}  
\|u_B\|_{L^2(A)}= |A|^{1/2}|u_B| = \|u_B\pm u\|_{L^2(A)}     \leq \|u-u_B\|_{L^2(B)} + \|u\|_{L^2(A)}. \end{equation}
 On $A$, $ |u(x)| \leq 2^s |x|^s|u(x)|$.  It follows that \begin{equation}\label{eq:2} \|u\|_{L^2(A)} \leq 2^s\||x|^su\|_{L^2(A)} \leq 2^s\||x|^su\|_{L^2(\mathbb R^d)}. \end{equation} Using the Poincar\'e inequality on the first term in the right-hand side of \eqref{eq:1} and \eqref{eq:2} for the second term, we therefore obtain \begin{equation}\label{eq:3} |u_B| \leq C_{d,s,|A|}  \left( \|u\|_{\dot H^s(B)} + \||x|^su\|_{L^2(\mathbb R^d)} \right ). \end{equation} We can now estimate the $L^2$-norm of $u$ inside the ball: \begin{equation}\label{eq:4} \begin{aligned} \|u\|_{L^2(B)} &\leq \|u-u_B\|_{L^2(B)} + \|u_B\|_{L^2(B)} \\ &= \|u-u_B\|_{L^2(B)} + |B|^{1/2}|u_B| \\ &\leq C_{d,s,|B|} \left( \|u\|_{\dot H^s(B)} + \||x|^su\|_{L^2(\mathbb R^d)} \right)\end{aligned} \end{equation}
 where we used again Poincar\'e  and \eqref{eq:3}. It remains to estimate $u$ outside $B$. 
 We have \begin{equation}\label{eq:5}  \|u\|_{L^2( B^c)}^2 = \int_{|x|\geq1}|u(x)|^2\,dx \leq \int_{|x|\geq1}|x|^{2s}|u(x)|^2\,dx \leq \||x|^su\|_{L^2(\mathbb R^d)}^2. \end{equation} Combining the estimates inside the unit ball, see \eqref{eq:4}, and outside the unit ball, see \eqref{eq:5},  yields \[ \|u\|_{L^2(\mathbb R^d)} \leq C_{d,s} \left( \|u\|_{\dot H^s(B)} + \||x|^su\|_{L^2(\mathbb R^d)} \right) \]
 and \eqref{frac_unc} is proved. The estimate \eqref{frac_unc-2} follows by $\|u\|_{\dot H^s(B)} \lesssim \|u\|_{\dot H^s(\mathbb R^d)}$ and a scaling argument. Indeed, for $u_\lambda=\lambda^{d/2}u(\lambda x)$, we have 
 \[
 \begin{aligned}
 \|u\|_{L^2(\mathbb R^d)}^2&=\|u_\lambda\|_{L^2(\mathbb R^d)}^2 \lesssim \left( \|u_\lambda\|_{\dot H^s(\mathbb R^d)}^2 + \||x|^s u_\lambda\|_{L^2(\mathbb R^d)}^2 \right)\\
 &\lesssim \left( \lambda^{2s}\|u \|_{\dot H^s(\mathbb R^d)}^2 + \lambda^{-2s}\||x|^s u \|_{L^2(\mathbb R^d)}^2 \right):=f(\lambda).
 \end{aligned}
 \]
  By optimizing in $\lambda$ we obtain that $\min f=f(\lambda_{\min})$, where $\lambda_{\min}=\|u \|_{\dot H^s(\mathbb R^d)}^{-1/(2s)}\||x|^s u \|_{L^2(\mathbb R^d)}^{1/(2s)}$
  and so
 $ \|u\|_{2}^2\lesssim \|u\|_{\dot H^s}\||x|^su\|_{L^2},$
which is \eqref{frac_unc-2}. 
\end{proof}

\section{Proof of Theorems \ref{thm:char} and \ref{thm:char-frac}}
We can now give the proof of the characterization Theorems \ref{thm:char} and \ref{thm:char-frac}. In order to prove scattering of a solution $\psi(t)$ to \eqref{NLS0}, we use Proposition \ref{prop:equiv-sca} and we work with the lens transform $\varphi$. Since $\Sigma^s$ is complete, proving \eqref{def-sca-lens} reduces to proving that the sequence $e^{-i\tau\mathcal H}\varphi(\tau)$ is a Cauchy family in $\Sigma^s$ as $\tau\to(\pi/4)^-$. 
We first prove the sufficient condition. Assume that $\sup_{t>0}\||J(t)|^s\psi(t)\|_{2}<\infty$. By Proposition \ref{prop:control-sigma-frac}, \begin{equation}\label{eq:varphi-sigma-s-bound} \sup_{0\leq\tau<\pi/4} \|\varphi(\tau)\|_{\Sigma^s} <\infty. \end{equation} 
By the Duhamel formulation \eqref{eq-duha},   for $ 0<\tau_1<\tau_2<\frac{\pi}{4}$, we write 
\begin{equation}\label{eq:duha-interaction-s}   e^{i\tau_2\mathcal H}\varphi(\tau_2)-e^{i\tau_1\mathcal H}\varphi(\tau_1)=-i\sum_{j=1}^N\mu_j \int_{\tau_1}^{\tau_2} e^{i\sigma\mathcal H} (\cos(2\sigma))^{-\delta(p_j)} (|\varphi|^{p_j-1}\varphi)(\sigma) \,d\sigma, 
\end{equation} 
where \begin{equation}\label{eq:dp}\delta(p):=2-\frac{d}{2}(p-1).
\end{equation}  Consequently, in order to prove that $e^{i\tau\mathcal H}\varphi(\tau)$ converges in $\Sigma^s$ as $\tau\to(\pi/4)^-$, it is enough to prove that the right-hand side of \eqref{eq:duha-interaction-s} tends to zero in $\Sigma^s$ as $ \tau_1,\tau_2\to(\pi/4)^-$.  Therefore, we estimate the nonlinear term on the right-hand side of \eqref{eq:duha-interaction-s}. We divide the analysis in three cases. 
\subsection*{Case $d>2s$}
Assume that $d>2s$, and set \begin{equation}\label{eq:def-q} q:=\frac{2d}{d-2s}.\end{equation}  For a fixed power $p_j$, $j\in\{1, \dots, N\}$, let $\rho_j$ be defined by  \begin{equation}\label{eq:rho-def} \frac1{\rho_j}:=\frac12-\frac{p_j-1}{2q}. \end{equation} Then $ \frac1{\rho_j'} = 1-\frac1{\rho_j} =  \frac{p_j-1}{q}+\frac1{\rho_j}$.  Thus, by the fractional chain rule and the corresponding weighted estimate, \begin{equation}\label{eq:nonlinear-Wsr-bootstrap} \||\varphi|^{p_j-1}\varphi \|_{\mathcal W^{s,\rho_j'}} \lesssim \|\varphi\|_{L^{q}}^{p_j-1} \|\varphi\|_{\mathcal W^{s,\rho_j}}. \end{equation} 
We now choose $\gamma_j$ so that $(\gamma_j,\rho_j)$ is a Schr\"odinger-admissible pair: therefore  \begin{equation}\label{eq:gamma-def} \gamma_j=\frac{8}{(p_j-1)(d-2s)}. \end{equation} 
Let $k_j$ be determined by  $  \frac1{\gamma_j'}=\frac1{k_j}+\frac1{\gamma_j}$.  It follows that \begin{equation}\label{def:k-time} \frac1{k_j} = 1-\frac2{\gamma_j} = 1-\frac{(p_j-1)(d-2s)}4  \iff k_j=\frac{4}{4-(p_j-1)(d-2s)}. \end{equation} Let $ \delta(p_j)$ as defined in \eqref{eq:dp}.  Consider the solution $\varphi$  to \eqref{NLS-lens}. For an interval $I=(\tau_0,\tau_1)\subset(0,\pi/4)$, the Strichartz estimate and \eqref{eq:nonlinear-Wsr-bootstrap} give 
\begin{equation*}\label{eq:bootstrap-Strichartz} 
\begin{aligned} 
\|\varphi\|_{L_\tau^{\gamma_j}(I;\mathcal W_x^{s,\rho_j})} &\lesssim \|\varphi(\tau_0)\|_{\Sigma^s} + \sum_{\ell=1}^N\| (\cos(2\cdot))^{-\delta(p_\ell)}|\varphi|^{p_\ell-1}\varphi \|_{L_\tau^{\gamma_\ell'}(I;\mathcal W_x^{s,\rho_\ell'})} \\ &\lesssim \|\varphi(\tau_0)\|_{\Sigma^s} + \sum_{\ell=1}^N \|(\cos(2\cdot))^{-\delta(p_\ell)}\|_{L_\tau^{k_\ell}(I)} \||\varphi|^{p_\ell-1}\varphi \|_{L_\tau^{\gamma_\ell}(I;\mathcal W_x^{s,\rho_\ell'})} \\ &\lesssim \|\varphi(\tau_0)\|_{\Sigma^s} + \sum_{\ell=1}^N \|(\cos(2\cdot))^{-\delta(p_\ell)}\|_{L_\tau^{k_\ell}(I)} \|\varphi\|_{L_\tau^\infty(I;L_x^{q})}^{p_\ell-1} \|\varphi\|_{L_\tau^{\gamma_\ell}(I;\mathcal W_x^{s,\rho_\ell})}. \end{aligned}
\end{equation*} 
Note that the  coefficients $\mu_j$ are absorbed into the implicit constants. By the Sobolev embedding $ H^s(\mathbb R^d)\hookrightarrow L^{q}(\mathbb R^d) $ and the uniform bound \eqref{eq:varphi-sigma-s-bound} we write  $\sup_{\tau\in(0,\pi/4)} \|\varphi(\tau)\|_{\Sigma^s} \leq C_0$, and then we have $ \|\varphi\|_{L_\tau^\infty(I;L_x^{q})} \lesssim C_0$. Consequently, 
\begin{equation}\label{eq:bootstrap-inequality-final} 
\|\varphi\|_{L_\tau^{\gamma_j}(I;\mathcal W_x^{s,\rho_j})} \lesssim C_0 +  \sum_{\ell=1}^NC_0^{p_\ell-1} \|(\cos(2\cdot))^{-\delta(p_\ell)}\|_{L_\tau^{k_\ell}(I)} \|\varphi\|_{L_\tau^{\gamma_\ell}(I;\mathcal W_x^{s,\rho_\ell})}. 
\end{equation} 
The condition for the time weight to belong to $L^{k_j}_\tau$   is 
$ \delta(p_j)k_j<1$.
Using \eqref{eq:dp} and \eqref{def:k-time}, this becomes $\left(2-\frac d2(p_j-1)\right) \frac4{4-(p_j-1)(d-2s)} <1. $ A direct calculation gives $ p_j>1+\frac4{d+2s}$. 
Therefore, for each $\ell\in\{1,\dots, N\}$, $\lim_{\tau_0\to(\pi/4)^-} \|(\cos(2\cdot))^{-\delta(p_\ell)}\|_{L_\tau^{k_\ell}((\tau_0,\pi/4))} =0$. Since we have finitely many nonlinear terms, we may thus choose $\tau_0 <\pi/4$, independent of $\ell$, sufficiently close to $\pi/4$ so that
\begin{equation}\label{eq-cost-abs} C \sum_{\ell=1}^N C_0^{p_\ell-1} \|(\cos(2\cdot))^{-\delta(p_\ell)}\|_{L_\tau^{k_\ell}((\tau_0,\pi/4))} \leq\frac1{2N}.\end{equation} 
Thus, summing \eqref{eq:bootstrap-inequality-final} over $j\in\{1,\dots,N\}$, we obtain \begin{align*} \sum_{j=1}^N \|\varphi\|_{L_\tau^{\gamma_j} (I;\mathcal W_x^{s,\rho_j})} \lesssim N C_0 + N\sum_{\ell=1}^N C_0^{p_\ell-1} \|(\cos(2\cdot))^{-\delta(p_\ell)}\|_{L_\tau^{k_\ell}(I)} \|\varphi\|_{L_\tau^{\gamma_\ell} (I;\mathcal W_x^{s,\rho_\ell})}, \end{align*} 
and we can absorb by using \eqref{eq-cost-abs} to  obtain \begin{equation}\label{eq:terminal-Strichartz-bound}\sum_{j=1}^N \|\varphi\|_ {L_\tau^{\gamma_j}((\tau_0,\pi/4);\mathcal W_x^{s,\rho_j})} \lesssim C_0. \end{equation} 
Finally, for $I=(\tau_1,\tau_2)\subset(\tau_0,\pi/4)$, the inhomogeneous Strichartz estimate and   \eqref{eq:terminal-Strichartz-bound} give 
\begin{equation} \label{eq:Duhamel-Cauchy-bootstrap}
\begin{aligned} 
& \|\sum_{j=1}^N \int_{\tau_1}^{\tau_2} e^{i\sigma\mathcal H} (\cos(2\sigma))^{-\delta(p_j)}  (|\varphi|^{p_j-1}\varphi)(\sigma) \,d\sigma \|_{\Sigma^s} \\ &\qquad\lesssim \sum_{j=1}^N\| (\cos(2\cdot))^{-\delta(p_j)}|\varphi|^{p_j-1}\varphi \|_ {L_\tau^{\gamma_j'}(I;\mathcal W_x^{s,\rho_j'})} \\ &\qquad\lesssim\sum_{j=1}^N \|(\cos(2\cdot))^{-\delta(p_j)}\|_{L_\tau^{k_j}(I)} \|\varphi\|_{L_\tau^\infty(I;L_x^{q})}^{p_j-1} \|\varphi\|_{L_\tau^{\gamma_j}(I;\mathcal W_x^{s,\rho_j})} \\ &\qquad\lesssim\sum_{j=1}^N C_0^{p_j-1} \|(\cos(2\cdot))^{-\delta(p_j)}\|_{L_\tau^{k_j}(I)} \|\varphi\|_ {L_\tau^{\gamma_j}((\tau_0,\pi/4);\mathcal W_x^{s,\rho_j})}
\\ &\qquad\lesssim\sum_{j=1}^N   \|(\cos(2\cdot))^{-\delta(p_j)}\|_{L_\tau^{k_j}(I)} . 
\end{aligned} 
\end{equation}
Since $ \|(\cos(2\cdot))^{-\delta(p_j)}\|_{L_\tau^{k_j}(\tau_1,\tau_2)} \to0 $ as $ \tau_1,\tau_2\to(\pi/4)^- $ and the sum is finite,  the right-hand side of \eqref{eq:Duhamel-Cauchy-bootstrap} tends to zero.

\subsection*{Case $d\leq2s$} Assume that $d\leq2s$. In this case, $ H^s \hookrightarrow L^q  \qquad\text{for every }2\leq q<\infty$. Notice that we use only finite Lebesgue exponents and do not rely on the algebra property of $H^s$ when $s>d/2$. For every $j\in\{1,\dots,N\}$, and for a fixed $q_j$ which will be chosen large enough, and define $\rho_j$ by 
\begin{equation}\label{eq:rho-d-leq-2s} \frac1{\rho_j} := \frac12-\frac{p_j-1}{2q_j} \iff  \frac1{\rho_j'} = \frac{p_j-1}{q_j}+\frac1{\rho_j}. 
\end{equation} 
By the fractional chain rule and the corresponding weighted estimate, valid for $0<s\leq1$, we have \begin{equation}\label{eq:nonlinear-d-leq-2s} \||\varphi|^{p_j-1}\varphi \|_{\mathcal W^{s,\rho_j'}} \lesssim \|\varphi\|_{L^{q_j}}^{p_j-1} \|\varphi\|_{\mathcal W^{s,\rho_j}}. \end{equation} Let $\gamma_j$ be such that $(\gamma_j,\rho_j)$  are related through the admissibility condition \eqref{eq:admissibility}.  In view of \eqref{eq:rho-d-leq-2s}, we have \begin{equation}\label{eq:gamma-d-leq-2s} \frac2{\gamma_j} = \frac{d(p_j-1)}{2q_j}. \end{equation} We choose $q_j$ sufficiently large so that $\gamma_j\geq2$. 
Define $k_j$ by  \begin{equation}\label{eq:k-d-leq-2s}\frac1{\gamma_j'} = \frac1{k_j}+\frac1{\gamma_j}\iff  \frac1{k_j}   = 1-\frac{d(p_j-1)}{2q_j}. \end{equation} 
Since $ p_j>1+\frac2d$, we have $ \delta(p_j) = 2-\frac d2(p_j-1) <1$. We claim that $q_j$ can be chosen so that \begin{equation}\label{eq:weight-integrability-d-leq-2s} \delta(p_j)k_j<1. \end{equation}By \eqref{eq:k-d-leq-2s}, condition \eqref{eq:weight-integrability-d-leq-2s} is equivalent to $1-\frac1{q_j}>\frac{2}{d(p_j-1)}$, which holds whenever $q_j$ is sufficiently large. Consequently, $ (\cos(2\tau))^{-\delta(p_j)} \in L^{k_j}(0,\pi/4)$. 
Consider an interval $I=(\tau_0,\tau_1)\subset(0,\pi/4)$. By the Strichartz estimates and \eqref{eq:nonlinear-d-leq-2s}, we obtain, for every $j\in\{1,\dots,N\}$, 
\begin{equation}\label{eq:bootstrap-d-leq-2s} 
\begin{aligned} \|\varphi\|_ {L_\tau^{\gamma_j}(I;\mathcal W_x^{s,\rho_j})} &\lesssim \|\varphi(\tau_0)\|_{\Sigma^s} + \sum_{\ell=1}^N \bigl\| (\cos(2\cdot))^{-\delta(p_\ell)} |\varphi|^{p_\ell-1}\varphi \bigr\|_ {L_\tau^{\gamma_\ell'} (I;\mathcal W_x^{s,\rho_\ell'})} \\ 
&\lesssim C_0+ \sum_{\ell=1}^N \|(\cos(2\cdot))^{-\delta(p_\ell)}\|_ {L_\tau^{k_\ell}(I)} \|\varphi\|_ {L_\tau^\infty(I;L_x^{q_\ell})}^{p_\ell-1} \|\varphi\|_ {L_\tau^{\gamma_\ell} (I;\mathcal W_x^{s,\rho_\ell})}\\ &\lesssim C_0+ \sum_{\ell=1}^N C_0^{p_\ell-1} \|(\cos(2\cdot))^{-\delta(p_\ell)}\|_ {L_\tau^{k_\ell}(I)} \|\varphi\|_ {L_\tau^{\gamma_\ell} (I;\mathcal W_x^{s,\rho_\ell})}. 
\end{aligned}\end{equation} 
Here we used \eqref{eq:varphi-sigma-s-bound} and the finite-exponent Sobolev embeddings $ H^s(\mathbb R^d)\hookrightarrow L^{q_\ell}$. Since the number of nonlinear terms is finite and $(\cos(2\tau))^{-\delta(p_\ell)} \in L^{k_\ell}(0,\pi/4),$ we may choose $\tau_0<\pi/4$, independently of $\ell$, sufficiently close to $\pi/4$ so that \[ C\sum_{\ell=1}^N C_0^{p_\ell-1} \|(\cos(2\cdot))^{-\delta(p_\ell)}\|_ {L_\tau^{k_\ell}(\tau_0,\pi/4)} \leq\frac1{2N}. \] Summing \eqref{eq:bootstrap-d-leq-2s} over $j\in\{1,\dots,N\}$ and absorbing the resulting term, we obtain \begin{equation}\label{eq:terminal-bound-d-leq-2s} \sum_{j=1}^N \|\varphi\|_ {L_\tau^{\gamma_j} ((\tau_0,\pi/4);\mathcal W_x^{s,\rho_j})} \lesssim C_0. \end{equation} 
Finally, let $ I=(\tau_1,\tau_2)\subset(\tau_0,\pi/4)$. By the inhomogeneous Strichartz estimate, \eqref{eq:nonlinear-d-leq-2s} and \eqref{eq:terminal-bound-d-leq-2s}, we have 
\begin{equation}\label{eq:Cauchy-d-leq-2s} 
\begin{aligned} & \left\| \sum_{j=1}^N\mu_j \int_{\tau_1}^{\tau_2} e^{i\sigma\mathcal H} (\cos(2\sigma))^{-\delta(p_j)} \bigl(|\varphi|^{p_j-1}\varphi\bigr)(\sigma) \,d\sigma \right\|_{\Sigma^s}  \\ &\qquad\lesssim \sum_{j=1}^N \|(\cos(2\cdot))^{-\delta(p_j)}\|_ {L_\tau^{k_j}(I)} \|\varphi\|_ {L_\tau^\infty(I;L_x^{q_j})}^{p_j-1} \|\varphi\|_ {L_\tau^{\gamma_j} (I;\mathcal W_x^{s,\rho_j})}  \\ & \qquad \lesssim \sum_{j=1}^N C_0^{p_j-1} \|(\cos(2\cdot))^{-\delta(p_j)}\|_ {L_\tau^{k_j}(I)} \|\varphi\|_ {L_\tau^{\gamma_j} ((\tau_0,\pi/4);\mathcal W_x^{s,\rho_j})}  \\ &\qquad\lesssim \sum_{j=1}^N \|(\cos(2\cdot))^{-\delta(p_j)}\|_ {L_\tau^{k_j}(I)}. 
\end{aligned}
\end{equation} 
For every $j$, the absolute continuity of the $L^{k_j}$ norm gives $\|(\cos(2\cdot))^{-\delta(p_j)}\|_ {L_\tau^{k_j}(\tau_1,\tau_2)} \to0 $ as $\tau_1,\tau_2\to\left(\frac\pi4\right)^-$.  Since the sum is finite, the right-hand side of \eqref{eq:Cauchy-d-leq-2s} tends to zero. Hence every Duhamel term tends to zero in $\Sigma^s$.

\subsection{Conclusion of the proof} We have proved that $e^{i\tau\mathcal H}\varphi(\tau)$ is a Cauchy family in $\Sigma^s$, which is  complete. Then, there exists $\overline{\varphi}^+\in\Sigma^s$ such that \[ \lim_{\tau\to(\pi/4)^-} \|e^{i\tau\mathcal H}\varphi(\tau)-\overline{\varphi}^+\|_{\Sigma^s}=0. \]
By the triangle inequality and the continuity of the harmonic linear propagator,  
 \[ \lim_{\tau\to(\pi/4)^-} \|\varphi(\tau)-\varphi^+ \|_{\Sigma^s} =0, \qquad \varphi^+:=e^{-i\frac{\pi}{4}\mathcal H}\overline{\varphi}^+, \]
 and by Proposition \ref{prop:equiv-sca}, the proof  of the sufficient implication is concluded. 
\medskip

\noindent We finally prove the converse implication. Suppose  that $\psi$ scatters forward in $\Sigma^s$. Then there exists $\psi^+\in\Sigma^s$ such that \[ \lim_{t\to\infty} \|U(-t)\psi(t)-\psi^+\|_{\Sigma^s}=0. \] By recalling that $|J(t)|^s=U(t)|x|^sU(-t)$, and by using the unitarity of $U(t)$ on $L^2$, we obtain $ \||J(t)|^s\psi(t)\|_2 = \||x|^sU(-t)\psi(t)\|_2$. Therefore, \[  | \||J(t)|^s\psi(t)\|_2 - \||x|^s\psi^+\|_2  |  \leq \| |x|^s(U(-t)\psi(t)-\psi^+) \|_2  \leq \|U(-t)\psi(t)-\psi^+\|_{\Sigma^s}.  \] It follows that \[ \lim_{t\to\infty} \||J(t)|^s\psi(t)\|_2 = \||x|^s\psi^+\|_2. \] In particular, there exists $t_0>0$ sufficiently large such that  \[ \sup_{t\geq t_0} \||J(t)|^s\psi(t)\|_2 <\infty. \] On $[0,t_0]$ the boundedness of $\||J(t)|^s\psi(t)\|_2 $ follows by continuity of the flow.


\section{Small-data scattering in $\Sigma^s$}

We now turn to the fractional small-data theory. As in the proof of Theorems \ref{thm:char} and \ref{thm:char-frac}, the key point is to obtain a uniform control of the solution $\varphi$ to \eqref{NLS-lens} in \[ L^\infty_\tau((0,\pi/4);\Sigma_x^s). \] The conclusion is a consequence of the following bootstrap lemma and the nonlinear estimates proved in the previous section.
\begin{lemma}\label{lem:bootstrap} Let $X:[0,\pi/4)\to[0,\infty)$ be continuous, and let $\eta\geq0$. Assume that \[ X(0)\leq C_0\eta \] and \[ X(T) \leq C_0\eta+C_1\sum_{j=1}^N X(T)^{p_j}, \qquad 0<T<\frac{\pi}{4}, \] where $p_j>1$. Then there exists \[ \varepsilon = \varepsilon(C_0,C_1,p_1,\dots,p_N,N)>0 \] such that, if $\eta\leq\varepsilon$, then \[ \sup_{0<T<\pi/4}X(T)\leq2C_0\eta. \] \end{lemma}
\begin{proof}
The proof is a standard continuity argument.
\end{proof}
We  now prove Theorem \ref{thm:small-data-frac}. 
\begin{proof}[Proof of Theorem \ref{thm:small-data-frac}] 
Let $\varphi$ be the lens transform of $\psi$, and let $[0,T_{\max})$ be its maximal forward interval of existence. We prove that, if $\|\varphi_0\|_{\Sigma^s}$ is sufficiently small, then \begin{equation}\label{eq:small-data-uniform-sigma} \sup_{0\leq\tau<T_{\max}} \|\varphi(\tau)\|_{\Sigma^s} \lesssim \|\varphi_0\|_{\Sigma^s}. \end{equation} This uniform bound, together with the continuation criterion, implies that $T_{\max}=\pi/4$. Theorem \ref{thm:char-frac} then yields forward scattering of $\psi$ in $\Sigma^s$.

\subsection*{Case $d>2s$} Set $q$ as in \eqref{eq:def-q}. For every $j\in\{1,\dots,N\}$, define $\rho_j$, $\gamma_j$, and $k_j$ as in \eqref{eq:rho-def}, \eqref{eq:gamma-def}, and \eqref{def:k-time}, respectively.  For $0<T<T_{\max}$, let 
\begin{equation*}\label{eq:X-small-data-supercritical} X(T) := \|\varphi\|_{L^\infty_\tau([0,T];\Sigma_x^s)} + \sum_{j=1}^N \|\varphi\|_ {L_\tau^{\gamma_j}([0,T];\mathcal W_x^{s,\rho_j})}. 
\end{equation*}
 By the Strichartz estimates, the fractional chain rule and the corresponding weighted estimate, exactly as in the proof of Theorem \ref{thm:char-frac}, we obtain \begin{equation}\label{eq:small-data-bootstrap-d-greater} X(T) \leq C\|\varphi_0\|_{\Sigma^s} + C\sum_{j=1}^N |\mu_j| \|(\cos(2\cdot))^{-\delta(p_j)}\|_{L^{k_j}_\tau(0,T)} X(T)^{p_j}, \end{equation} where the constant $C$ is independent of $T$. The assumption $ p_j>1+\frac4{d+2s} $ is equivalent to $ \delta(p_j)k_j<1. $ Consequently, $ (\cos(2\tau))^{-\delta(p_j)} \in L^{k_j}(0,\pi/4)$, and hence \[ \sup_{0<T<\pi/4} \|(\cos(2\cdot))^{-\delta(p_j)}\|_{L^{k_j}(0,T)} <\infty. \] It follows from \eqref{eq:small-data-bootstrap-d-greater} that \[ X(T) \leq C_0\|\varphi_0\|_{\Sigma^s} + C_1\sum_{j=1}^N X(T)^{p_j}, \] with constants $C_0,C_1$ independent of $T$. Lemma \ref{lem:bootstrap} therefore gives, provided $\|\varphi_0\|_{\Sigma^s}$ is sufficiently small, \begin{equation*}\label{eq:small-data-X-bound-supercritical} \sup_{0<T<\pi/4}X(T) \lesssim \|\varphi_0\|_{\Sigma^s}. \end{equation*} 
 In particular, \eqref{eq:small-data-uniform-sigma} holds. 
\subsection*{Case $d\leq2s$} For every $j$, choose a sufficiently large finite exponent $q_j$ and define $\rho_j$, $\gamma_j$, and $k_j$ as in \eqref{eq:rho-d-leq-2s}, \eqref{eq:gamma-d-leq-2s}, and \eqref{eq:k-d-leq-2s}, respectively.   Since $ H^s(\mathbb R^d)\hookrightarrow L^{q_j}(\mathbb R^d) $ for every finite $q_j$, the same Strichartz and nonlinear estimates of the previous section give, with \[ X(T) := \|\varphi\|_{L^\infty_\tau([0,T];\Sigma_x^s)} + \sum_{j=1}^N \|\varphi\|_ {L_\tau^{\gamma_j}([0,T];\mathcal W_x^{s,\rho_j})}, \] the inequality \begin{equation}\label{eq:small-data-bootstrap-critical} X(T) \leq C\|\varphi_0\|_{\Sigma^s} + C\sum_{j=1}^N |\mu_j| \|(\cos(2\cdot))^{-\delta(p_j)}\|_{L^{k_j}_\tau(0,T)} X(T)^{p_j}. \end{equation} In this case, the assumption $ p_j>1+\frac2d $ gives $\delta(p_j)<1$. We can choose $q_j$ sufficiently large so that $ \delta(p_j)k_j<1$. Therefore, \[ \sup_{0<T<\pi/4} \|(\cos(2\cdot))^{-\delta(p_j)}\|_{L^{k_j}(0,T)} <\infty.
\]
Applying Lemma \ref{lem:bootstrap} to \eqref{eq:small-data-bootstrap-critical}, we again obtain \[ \sup_{0<T<\pi/4}X(T) \lesssim \|\varphi_0\|_{\Sigma^s}, \] provided $\|\varphi_0\|_{\Sigma^s}$ is sufficiently small. In particular, \eqref{eq:small-data-uniform-sigma} holds. 
  \end{proof}
 
 \subsection{Proof of Corollary \ref{coro-sd}}

We  now prove Corollary \ref{coro-sd}. It combines the use of the fractional uncertainty principle, see Proposition \ref{prop:u-p}  and the small data scattering result in Theorem \ref{thm:small-data-frac}.
\begin{proof}[Proof of Corollary \ref{coro-sd}]
Consider the equation-preserving scaling as defined in \eqref{NLS-s-scaling}, and  define $\Gamma_\lambda(f)(x)=\lambda^{2/(p-1)}f(\lambda x)$, hence $\psi_\lambda(t,x)=\Gamma_\lambda(\psi(\lambda^2t))$. By using  \eqref{frac_unc-2} we have
\[\begin{aligned}
\|\Gamma_\lambda(\psi_0)\|_{\Sigma^s}^2&=\|\psi_{0,\lambda}\|_2^2+\||\nabla|^s\psi_{0,\lambda}\|_2^2+\||\cdot |^s\psi_\lambda\|_2^2\\
&=\lambda^{\frac{4}{p-1}-d}\|\psi_0\|_2^2+\lambda^{2s-d+\frac{4}{p-1}}\||\nabla |^s\psi_0\|_2^2+\lambda^{-2s-d+\frac{4}{p-1}}\||\cdot |^s\psi_0\|_2^2\\
&\leq C\lambda^{\frac{4}{p-1}-d}\||x|^s\psi_0\|_2\||\nabla|^s\psi_0\|_2+\lambda^{2s-d+\frac{4}{p-1}}\||\nabla |^s\psi_0\|_2^2\\
&\qquad+\lambda^{-2s-d+\frac{4}{p-1}}\||\cdot |^s\psi_0\|_2^2\\
&:=CAB\lambda^{\theta} +\lambda^{\theta+2s}A^2+\lambda^{\theta-2s}B^2=f(\lambda), 
\end{aligned}
\]
where 
\[
A=\||\nabla|^s\psi_0\|_{2}, \quad  B=\||x|^s\psi_0\|_2, \quad \theta=\frac{4}{p-1}-d.
\]
We optimize as follows: 
 \[
\begin{aligned}
&f'(\lambda)=\lambda^{\theta-1-2s}\left(\theta \lambda^{2s} CAB+(\theta+2s)\lambda^{4s}A^2+(\theta-2s) B^2  \right)=0\\
 &\iff \theta CAB\lambda^{2s}+(\theta+2s)A^2\lambda^{4s}+(\theta-2s) B^2  )=0,\\
& \iff \lambda=\left(\frac{-\theta C\pm\sqrt{\theta^2(C^2-4)+16s^2}}{2\theta+4s}\frac{B}{A}\right)^{1/(2s)}
\end{aligned}
\]
Note that  under the hypothesis    $1+\max\left\{\frac2d,\frac4{d+2s}\right\} <p <1+\frac4d$, we have $ 0<\theta<2s. $
Moreover, $f(\lambda)\to+\infty$ for $\lambda\to0$ and $f(\lambda)\to+\infty $ for $\lambda\to+\infty.$ Thus, $f$ attains its minimum at the  positive root \[ \lambda_{\min} = \left( \frac{B}{A} \frac{ \sqrt{\theta^2C^2+4(4s^2-\theta^2)} -\theta C }{ 2(\theta+2s) }\right)^{1/(2s)}. \] Note that since $0<\theta<2s$, one has $\sqrt{\theta^2C^2+4(4s^2-\theta^2)} > \theta C$. Set $K=\frac{ \sqrt{\theta^2C^2+4(4s^2-\theta^2)} -\theta C }{ 2(\theta+2s) }$. Substituting $\lambda_{\min}$ into $f$, we obtain 
\begin{align*} f(\lambda_{\min})&= \left( K \frac BA \right)^{\frac{\theta}{2s}} \left( CAB +A^2K \frac BA +B^2K ^{-1}\frac AB \right) \\ &= \left( K \frac BA \right)^{\frac{\theta}{2s}} AB \left( C+K  +K ^{-1} \right) \\ &=K ^{\frac{\theta}{2s}} A^{1-\frac{\theta}{2s}} B^{1+\frac{\theta}{2s}}  \left( C+K  +K ^{-1} \right). \end{align*} Therefore, \[   \min_{\lambda>0}f(\lambda) = \tilde K_{\theta,s,C} A^{1-\frac{\theta}{2s}} B^{1+\frac{\theta}{2s}}= \tilde K_{\theta,s,C}\||\nabla|^s\psi_0\|_2^{1-\frac{\theta}{2s}}\||x|^s\psi_0\|_2^{1+\frac{\theta}{2s}}. \] 
This implies that 
\[
\|\Gamma_{\lambda_{\min}}(\psi_0)\|_{\Sigma^s}\leq\min_{\lambda>0}\|\Gamma_\lambda(\psi_0)\|_{\Sigma^s}\lesssim \||\nabla|^s\psi_0\|_2^{\frac12-\frac{\theta}{4s}}\||x|^s\psi_0\|_2^{\frac12+\frac{\theta}{4s}}
\]
and so, since by hypothesis $\||\nabla|^s\psi_0\|_2^{\frac12-\frac{\theta}{4s}}\||x|^s\psi_0\|_2^{\frac12+\frac{\theta}{4s}}$ is small, we conclude that  the solution to \eqref{NLS-s} with initial datum $\tilde\psi=\Gamma_{\lambda_{\min}}(\psi_0)$ scatters to a profile $\psi^+$. By the fact that $\Gamma_\lambda U(\lambda^2t)=U(t)\Gamma_\lambda$, we have that the solution to \eqref{NLS-s} with initial datum $\psi_0$ scatters. Indeed, 
\[
\begin{aligned}
\|U(-\lambda_{\min}^2t)\psi(\lambda_{\min}^2t)-\Gamma_{\lambda_{\min}^{-1}}\psi^+\|_{\Sigma^s}&=\|\Gamma_{\lambda_{\min}^{-1}}\Gamma_{\lambda_{\min}}U(-\lambda_{\min}^2t)\psi(\lambda_{\min}^2t) -\Gamma_{\lambda_{\min}^{-1}}\psi^+\|_{\Sigma^s}
\\&\lesssim_{\lambda_{\min}}\|\Gamma_{\lambda_{\min}}U(-\lambda_{\min}^2t)\psi(\lambda_{\min}^2t)-\psi^+\|_{\Sigma^s}\\
&=\|U(-t)\Gamma_{\lambda_{\min}}\psi(\lambda_{\min}^2t)-\psi^+\|_{\Sigma^s}\\
&=\|U(-t) \psi_{\lambda_{\min}}(t)-\psi^+\|_{\Sigma^s}. 
\end{aligned}
\]
Hence, by setting $\lambda^2t=t'$, we have 
\[\|U(-t')\psi(t')-\Gamma_{\lambda_{\min}^{-1}}\psi^+\|_{\Sigma^s}\to 0\]
as $t'\to\infty$, which is the scattering property for $\psi$ with scattering state $\Gamma_{\lambda_{\min}^{-1}}\psi^+$.
\end{proof}
 
\subsection*{Acknowledgments} This work is supported by the University of Pisa and the Ministry of University and Research under the ``PRO3 - Programmazione triennale 2024-2026''. LF, VG, and NV are members of the GNAMPA group of the INdAM. LF and VG   are partially supported by the INdAM-GNAMPA Project, cod. CUP E53C25002010001. 
V.G.  is partially supported    by the Institute of Mathematics and Informatics, Bulgarian Academy of Sciences and by Waseda University of Tokyo. 


\end{document}